\documentclass[12pt,twoside]{amsart}
\usepackage{amsfonts}
\usepackage{amssymb}
\usepackage{amsmath}
\usepackage{fancyhdr}

\theoremstyle{plain}
\newtheorem{thm}{Theorem}[section]
\newtheorem{theorem}[thm]{Theorem}
\newtheorem{lemma}[thm]{Lemma}
\newtheorem{corollary}[thm]{Corollary}
\newtheorem{proposition}[thm]{Proposition}

\theoremstyle{definition}

\newtheorem{remark}[thm]{Remark}

\newtheorem{defn-thm}[thm]{Definition-Theorem}

\numberwithin{equation}{section}
\catcode`\@=11
\def\opn#1#2{\def#1{\mathop{\kern0pt\fam0#2}\nolimits}}
\def\underrightarrow{\mathpalette\underrightarrow@}
\def\underrightarrow@#1#2{\vtop{\ialign{$##$\cr
 \hfil#1#2\hfil\cr\noalign{\nointerlineskip}%
 #1{-}\mkern-6mu\cleaders\hbox{$#1\mkern-2mu{-}\mkern-2mu$}\hfill
 \mkern-6mu{\to}\cr}}}

\def\underleftarrow{\mathpalette\underleftarrow@}
\def\underleftarrow@#1#2{\vtop{\ialign{$##$\cr
 \hfil#1#2\hfil\cr\noalign{\nointerlineskip}#1{\leftarrow}\mkern-6mu
 \cleaders\hbox{$#1\mkern-2mu{-}\mkern-2mu$}\hfill
 \mkern-6mu{-}\cr}}}
\let\amp@rs@nd@\relax
\newdimen\ex@
\newdimen\bigaw@
\newdimen\minaw@
\newdimen\minCDaw@
\newif\ifCD@
\def\minCDarrowwidth#1{\minCDaw@#1}

\def\@CD{\def\A##1A##2A{\llap{$\vcenter{\hbox
 {$\scriptstyle##1$}}$}\Big\uparrow\rlap{$\vcenter{\hbox{%
$\scriptstyle##2$}}$}&&}%
\def\V##1V##2V{\llap{$\vcenter{\hbox
 {$\scriptstyle##1$}}$}\Big\downarrow\rlap{$\vcenter{\hbox{%
$\scriptstyle##2$}}$}&&}%
\def\={&\hskip.5em\mathrel
 {\vbox{\hrule width\minCDaw@\vskip3\ex@\hrule width
 \minCDaw@}}\hskip.5em&}%
\def\verteq{\Big\Vert&&}%
\def\noarr{&&}%
\def\vspace##1{\noalign{\vskip##1\relax}}\relax\iffalse{%
\fi\let\amp@rs@nd@&\iffalse}\fi
 \CD@true\vcenter\bgroup\relax\iffalse{%
\fi\let\\=\cr\iffalse}\fi\tabskip\z@skip\baselineskip20\ex@
 &\hfill$\m@th##$\hfill\cr}
\def\@endCD{\cr\egroup\egroup}
\def\>#1>#2>{\amp@rs@nd@\setbox\z@\hbox{$\scriptstyle
 \;{#1}\;\;$}\setbox\@ne\hbox{$\scriptstyle\;{#2}\;\;$}\setbox\tw@
 \hbox{$#2$}\ifCD@
 \global\bigaw@\minCDaw@\else\global\bigaw@\minaw@\fi
 \ifdim\wd\z@>\bigaw@\global\bigaw@\wd\z@\fi
 \ifdim\wd\@ne>\bigaw@\global\bigaw@\wd\@ne\fi
 \ifCD@\hskip.5em\fi
 \ifdim\wd\tw@>\z@
 \mathrel{\mathop{\hbox to\bigaw@{\rightarrowfill}}\limits^{#1}_{#2}}\else
 \mathrel{\mathop{\hbox to\bigaw@{\rightarrowfill}}\limits^{#1}}\fi
 \ifCD@\hskip.5em\fi\amp@rs@nd@}
\def\<#1<#2<{\amp@rs@nd@\setbox\z@\hbox{$\scriptstyle
 \;\;{#1}\;$}\setbox\@ne\hbox{$\scriptstyle\;\;{#2}\;$}\setbox\tw@
 \hbox{$#2$}\ifCD@
 \global\bigaw@\minCDaw@\else\global\bigaw@\minaw@\fi
 \ifdim\wd\z@>\bigaw@\global\bigaw@\wd\z@\fi
 \ifdim\wd\@ne>\bigaw@\global\bigaw@\wd\@ne\fi
 \ifCD@\hskip.5em\fi
 \ifdim\wd\tw@>\z@
 \mathrel{\mathop{\hbox to\bigaw@{\leftarrowfill}}\limits^{#1}_{#2}}\else
 \mathrel{\mathop{\hbox to\bigaw@{\leftarrowfill}}\limits^{#1}}\fi
 \ifCD@\hskip.5em\fi\amp@rs@nd@}
\def\@CDS{\def\A##1A##2A{\llap{$\vcenter{\hbox
 {$\scriptstyle##1$}}$}\Big\uparrow\rlap{$\vcenter{\hbox{%
$\scriptstyle##2$}}$}&}%
\def\V##1V##2V{\llap{$\vcenter{\hbox
 {$\scriptstyle##1$}}$}\Big\downarrow\rlap{$\vcenter{\hbox{%
$\scriptstyle##2$}}$}&}%
\def\={&\hskip.5em\mathrel
 {\vbox{\hrule width\minCDaw@\vskip3\ex@\hrule width
 \minCDaw@}}\hskip.5em&}
\def\verteq{\Big\Vert&}
\def\novarr{&}
\def\noharr{&&}
\def\SE##1E##2E{\slantedarrow(0,18)(4,-3){##1}{##2}&}
\def\SW##1W##2W{\slantedarrow(24,18)(-4,-3){##1}{##2}&}
\def\NE##1E##2E{\slantedarrow(0,0)(4,3){##1}{##2}&}
\def\NW##1W##2W{\slantedarrow(24,0)(-4,3){##1}{##2}&}
\def\slantedarrow(##1)(##2)##3##4{%
\thinlines\unitlength1pt\lower 6.5pt\hbox{\begin{picture}(24,18)%
\put(##1){\vector(##2){24}}%
\put(0,8){$\scriptstyle##3$}%
\put(20,8){$\scriptstyle##4$}%
\end{picture}}}
\def\vspace##1{\noalign{\vskip##1\relax}}\relax\iffalse{%
\fi\let\amp@rs@nd@&\iffalse}\fi
 \CD@true\vcenter\bgroup\relax\iffalse{%
\fi\let\\=\cr\iffalse}\fi\tabskip\z@skip\baselineskip20\ex@
 &\hfill$\m@th##$\hfill\cr}
\def\@endCDS{\cr\egroup\egroup}
\newdimen\TriCDarrw@
\newif\ifTriV@

\def\@TriCDV{\TriV@true\def\TriCDpos@{6}\@TriCD}
\def\@TriCDA{\TriV@false\def\TriCDpos@{10}\@TriCD}
\def\@TriCD#1#2#3#4#5#6{%
\setbox0\hbox{$\ifTriV@#6\else#1\fi$} \TriCDarrw@=\wd0
\advance\TriCDarrw@ 24pt \advance\TriCDarrw@ -1em
\def\SE##1E##2E{\slantedarrow(0,18)(2,-3){##1}{##2}&}
\def\SW##1W##2W{\slantedarrow(12,18)(-2,-3){##1}{##2}&}
\def\NE##1E##2E{\slantedarrow(0,0)(2,3){##1}{##2}&}
\def\NW##1W##2W{\slantedarrow(12,0)(-2,3){##1}{##2}&}
\def\slantedarrow(##1)(##2)##3##4{\thinlines\unitlength1pt
\lower 6.5pt\hbox{\begin{picture}(12,18)%
\put(##1){\vector(##2){12}}%
\put(-4,\TriCDpos@){$\scriptstyle##3$}%
\put(12,\TriCDpos@){$\scriptstyle##4$}%
\end{picture}}}
\def\={\mathrel {\vbox{\hrule
   width\TriCDarrw@\vskip3\ex@\hrule width
   \TriCDarrw@}}}
\def\>##1>>{\setbox\z@\hbox{$\scriptstyle
 \;{##1}\;\;$}\global\bigaw@\TriCDarrw@
 \ifdim\wd\z@>\bigaw@\global\bigaw@\wd\z@\fi
 \hskip.5em
 \mathrel{\mathop{\hbox to \TriCDarrw@
{\rightarrowfill}}\limits^{##1}}
 \hskip.5em}
\def\<##1<<{\setbox\z@\hbox{$\scriptstyle
 \;{##1}\;\;$}\global\bigaw@\TriCDarrw@
 \ifdim\wd\z@>\bigaw@\global\bigaw@\wd\z@\fi
 \mathrel{\mathop{\hbox to\bigaw@{\leftarrowfill}}\limits^{##1}}
 }
 \CD@true\vcenter\bgroup\relax\iffalse{\fi\let\\=\cr\iffalse}\fi
 \tabskip\z@skip\baselineskip20\ex@
 \ifTriV@
 \halign\bgroup
 &\hfill$\m@th##$\hfill\cr
#1&\multispan3\hfill$#2$\hfill&#3\\
&#4&#5\\
&&#6\cr\egroup%
\else
 \halign\bgroup
 &\hfill$\m@th##$\hfill\cr
&&#1\\%
&#2&#3\\
#4&\multispan3\hfill$#5$\hfill&#6\cr\egroup \fi}
\def\@endTriCD{\egroup}

\newcounter{Myenumi}
{\begin{list}{}{\usecounter{Myenumi}%
\settowidth{\leftmargin}{2.n}\settowidth{\labelwidth}{2.n}%
\setlength{\labelsep}{0pt}}}{\end{list}}
\newcounter{Myenumii}
{\begin{list}{}{\usecounter{Myenumii}%
\settowidth{\leftmargin}{a)n}\settowidth{\labelwidth}{a)n}%
\setlength{\labelsep}{0pt}}}{\end{list}}
\newcounter{Myenumiii}
{\begin{list}{}{\usecounter{Myenumiii}%
\settowidth{\leftmargin}{iv.n}\settowidth{\labelwidth}{iv.n}%
\setlength{\labelsep}{0pt}}}{\end{list}}

{\end{list}}

{\begin{list}{}{%
\settowidth{\leftmargin}{2.n}\settowidth{\labelwidth}{2.n}%
\setlength{\labelsep}{0pt}}}{\end{list}}

\newsymbol\onto 1310

\begin{document}

\title{Regularity of Leray-Hopf solutions to Navier-Stokes equations}

\author{Jian Zhai \\ \tiny{Department of Mathematics, Zhejiang University, Hangzhou
310027, PRC} } \maketitle

\footnotetext{This work is supported by NSFC No.10571157. email:
jzhai@zju.edu.cn}

%\begin{center}
%\email{jzhai@zju.edu.cn}
%\end{center}

%\tableofcontents

%\newpage

\begin{abstract}
An upper bound of blow up rate for impressible Navier-Stokes
equations with small data in $L^2(\mathbb R^3)$ is obtained.
\end{abstract}

\section{Introduction}

We consider the blow up rate of weak solutions to impressible
Navier-Stokes equations
\begin{equation}\label{1.1}
\left\{\begin{aligned}
\partial_tu-\Delta u+u\cdot\nabla u+\nabla p=0,\quad\text{in}\quad \mathbb R^3\times(0,T)\\
\text{div}u=0,\quad\text{in}\quad \mathbb R^3\times(0,T)\\
u(x,0)=u_0(x),\quad\text{in}\quad \mathbb R^3
\end{aligned}\right.
\end{equation}
where  $u$ and $p$ denote the unknown velocity and pressure of
incompressible fluid respectively.\\

In this paper, we shall estimate the upper bound of blow up rate for
the
Navier-Stokes equations.\\

\begin{theorem}\label{thm1.1} There is $\delta>0$ such that if
$\|u_0\|_{L^2(\mathbb R^3)}\leq\delta$, and if $u$ is a Leray-Hopf
solution to the problem (\ref{1.1}) and blows up at $t=T$,  then for
any small $\epsilon>0$, there is $t_0\in(0,T)$, such that
\begin{equation}\label{1.2}
\|u(t)\|_{L^\infty(\mathbb R^3)}\leq
\frac{\epsilon}{(T-t)^{1/2}},\quad \text{for all}\quad
t\in(t_0,T).
\end{equation}

\end{theorem}

\vspace{0.5cm}

Here $u: \,\, (x,t)\in\mathbb R^3\times(0,T)\to\mathbb R^3$ is
called a weak solution of (\ref{1.1}) if it is a Leray-Hopf
solution. Precisely, it satisfies
\begin{eqnarray*}
&&(1)\quad u\in L^\infty(0,T;L^2(\mathbb R^3))\cap L^2(0,T;H^{1}(\mathbb R^3)),\\
&&(2)\quad  \text{div} u=0 \quad\text{in}\quad \mathbb R^3\times(0,T),\\
&&(3)\quad \int_0^T\int_{\mathbb
R^3}\{-u\cdot\partial_t\phi+\nabla u\cdot\nabla\phi+(u\cdot\nabla
u)\cdot\phi\}dxdt=0
\end{eqnarray*}
for all $\phi\in C^\infty_0(\mathbb R^3\times(0,T))$ with
div$\phi=0$ in $\mathbb R^3\times(0,T)$.

 Combining Theorem 1.1 with my former result in \cite{[Z]}, we have

\begin{corollary}\label{cor1.2} There is $\delta>0$ such that if
$\|u_0\|_{L^2(\mathbb R^3)}\leq\delta$, and if $u$ is a Leray-Hopf
solution of the Navier-Stokes equations (\ref{1.1}), then $u$ is
regular in $\mathbb R^3\times(0,\infty)$.
\end{corollary}

 \vspace{0.5cm}

Since Leray(1934)\cite{[L]} and Hopf(1951)\cite{[H]} obtained the
global existence of weak solutions, it has been a fundamental open
problem to prove the uniqueness and regularity of weak solutions
to the Navier-Stokes equations.\\

\section{ Energy estimates}

As in \cite{[GK]}\cite{[GK1]}\cite{[GK2]} where Giga and Kohn
introduced similar transformations for the blow-up problem of
semi-linear heat equations, we apply
\begin{equation}\label{2.1}
y=\frac1{(T-t)^{1/2}}x,\quad \tau=-\ln (T-t),\quad
w(y,\tau)=(T-t)^{1/2}u(x,t),
\end{equation}
 to (\ref{1.1}) and consider the following new problem
\begin{equation}\label{2.2}
\left\{\begin{aligned}
\partial_\tau
w=\Delta_yw-\frac{y}2\cdot\nabla_yw-\frac12
w-w\cdot\nabla_yw-\nabla_yq,\quad \forall y\in\mathbb R^3,\quad
\tau>-\ln T\\
\text{div}_y\,\, w(y,\tau)=0,\quad\text{in}\quad \mathbb R^3\times(-\ln T,\infty) \\
w(y,-\ln T)=T^{1/2}u_0(T^{1/2}y),\quad\text{in}\quad\mathbb R^3
\end{aligned}\right.
\end{equation}
where
$$
q(y,\tau)=(T-t)p(x,t).
$$

Without loss generality, in this section we take $T=1$.
Multiplying the first one of (\ref{2.2}) by $w$ and integrating it
over $\mathbb R^3$, by using the second equation of (\ref{2.2}) we
have
\begin{equation}\label{2.3}
\begin{aligned}
\frac12\int_{\mathbb R^3}\partial_\tau |w(y,\tau)|^2dy
=(-1)\int_{\mathbb
R^3}|\nabla_yw(y,\tau)|^2-\frac14|w(y,\tau)|^2dy\\
-\frac14\int_{\mathbb R^3}\text{div}\,\, (y|w(y,\tau)|^2)dy .
\end{aligned}
\end{equation}
Noting that
\begin{equation}\label{2.4}
\int_{\mathbb R^3}\text{div}\,\,
(y|w(y,\tau)|^2)dy=\lim_{R\to\infty}\int_{\partial
B_R}|y||w(y,\tau)|^2d\sigma(y)\geq 0
\end{equation}
we obtain

\begin{lemma}\label{lem2.1} For any $\tau>0$, we have
\begin{equation}\label{2.5}
\frac12\frac{d}{d\tau}\|w(\tau)\|_{L^2(\mathbb R^3)}^2\leq
(-1)\{\|\nabla_yw(\tau)\|_{L^2(\mathbb
R^3)}^2-\frac14\|w(\tau)\|_{L^2(\mathbb R^3)}^2\}.
\end{equation}

\end{lemma}

\vspace{0.5cm}

Furthermore, we take differential in the equations of (\ref{2.2})
and obtain
\begin{equation}\label{2.6}
\begin{aligned}
\partial_\tau\partial_jw
=\Delta\partial_jw-\frac12y\cdot\nabla\partial_jw-\partial_jw\\
-(\partial_jw\cdot\nabla)w-(w\cdot\nabla)\partial_jw-\nabla_y\partial_jq.
\end{aligned}
\end{equation}
By the same strategy as in the proof of Lemma \ref{lem2.1}, from
(\ref{2.6}) as well as the equation
\begin{eqnarray*}
\partial_\tau \Delta w=\Delta^2 w-\frac12 (y\cdot\nabla) \Delta
w-\frac32\Delta w-\Delta((w\cdot\nabla)w)-\nabla\Delta q
\end{eqnarray*}
by taking twice differential in (\ref{2.2}), we have

\begin{lemma}\label{lem2.2} For all $\tau>0$
\begin{equation}\label{2.7}
\begin{aligned}
\frac{d}{d\tau}\int_{\mathbb R^3}|\nabla w(y,\tau)|^2dy \leq
-2\int_{\mathbb R^3}|\nabla^2w(y,\tau)|^2dy-\frac12\int_{\mathbb
R^3}|\nabla w(y,\tau)|^2dy\\
-2\sum_{j,k,l=1}^3\int_{\mathbb
R^3}\partial_jw_k(y,\tau)\partial_jw_l(y,\tau)\partial_lw_k(y,\tau)dy
\end{aligned}
\end{equation}
and
\begin{equation}\label{2.8}
\begin{aligned}
\frac{d}{d\tau}\int_{\mathbb R^3}|\Delta w(y,\tau)|^2dy \leq
-2\int_{\mathbb R^3}|\nabla\Delta
w(y,\tau)|^2dy-\frac32\int_{\mathbb
R^3}|\Delta w(y,\tau)|^2dy\\
-2\int_{\mathbb R^3}(\Delta
w(y,\tau))\cdot\Delta((w(y,\tau)\cdot\nabla)w(y,\tau))dy.
\end{aligned}
\end{equation}

\end{lemma}

\vspace{0.5cm}

\begin{remark}\label{rmk2.3} (1) For any $t_1>0$, there is
$t_0\in(0,t_1)$ such that $u(\cdot,t_0)\in H^1(\mathbb R^3)$. With
the initial data $u(x,t_0)$, the Leray-Hopf solution $u(x,t)$ is
regular at least in a short time interval after $t_0$ (see
\cite{[L]}\cite{[S]}). We are discussing the blow-up problem for
these short time regular solutions.

(2) As a blow-up argument, we assume that $u(x,t)$ is bounded for
$t<T$ and blows up at $t=T$. As a direct corollary, we can prove
that $\|u(t)\|_{H^3(\mathbb R^3)}$ and
$\partial_t\|u(t)\|_{H^m(\mathbb R^3)}$ ($m=0,1,2$), as well as
$\|\partial_tu(t)\|_{L^2(\mathbb R^3)},\,\,
\|\partial_t\nabla_xu(t)\|_{L^2(\mathbb R^3)}$ are bounded for
$t<T$. So we have the same results for $\|w(\tau)\|_{H^3(\mathbb
R^3)}$ and $\partial_\tau\|w(\tau)\|_{H^m(\mathbb R^3)}$ ($m=0,1,2$)
for $\tau<\infty$, as well as the similar results for $q$ by the
boundedness of Riesz transformation.

(3) Since $u(x,t),\,\,\,\partial_tu(x,t)\in L^2(\mathbb R^3)$ for
$t<T$,
$$
\int_0^t\int_{\mathbb R^3}|\partial_hu(x,h)||u(x,h)|dxdh<\infty,
$$
we can use Fubini theorem to obtain
\begin{equation}\label{2.9}
\begin{aligned}
2\int_{\mathbb R^3}\partial_tu(x,t)\cdot
u(x,t)dx=\frac{d}{dt}\int_0^t\int_{\mathbb
R^3}\partial_h|u(x,h)|^2dxdh\\
=\frac{d}{dt}\int_{\mathbb
R^3}\int_0^t\partial_h|u(x,h)|^2dhdx=\frac{d}{dt}\int_{\mathbb
R^3}|u(x,t)|^2dx.
\end{aligned}\end{equation}
Noting that
$$
\partial_tu(x,t)=(T-t)^{-\frac32}\{\partial_\tau
w(\frac{x}{(T-t)^{1/2}},\tau)+\frac{x}{2(T-t)^{1/2}}\cdot\nabla_yw(\frac{x}{(T-t)^{1/2}},\tau)+\frac12w(\frac{x}{(T-t)^{1/2}},\tau)\}
$$
where $\tau=(-)\ln(T-t)$, from
\begin{equation}\label{2.10}
\int_{\mathbb
R^3}|\partial_tu(x,t)|^2dx=(T-t)^{-\frac32}\int_{\mathbb
R^3}|\partial_\tau w(y,\tau)+\frac{y}2\cdot\nabla_y
w(y,\tau)+\frac12w(y,\tau)|^2dy
\end{equation}
and
\begin{equation}\label{2.11}
\int_{\mathbb R^3}|u(x,t)|^2dx=(T-t)^{1/2}\int_{\mathbb
R^3}|w(y,\tau)|^2dy
\end{equation}
we have for $t<T$
\begin{equation}\label{2.12}
\int_{\mathbb R^3}|\partial_\tau
w(y,\tau)+\frac{y}2\cdot\nabla_yw(y,\tau)|^2dy<\infty.
\end{equation}
Moreover, from (\ref{2.9}), we get
\begin{equation}\label{2.13}
\begin{aligned}
(T-t)^{-\frac12}\{\partial_\tau\int_{\mathbb
R^3}|w(y,\tau)|^2dy-\frac12\int_{\mathbb R^3}|w(y,\tau)|^2dy\}\\
=\frac{d}{dt}\{(T-t)^{\frac12}\int_{\mathbb
R^3}|w(y,\tau)|^2dy\}=\frac{d}{dt}\int_{\mathbb
R^3}|u(x,t)|^2dx=2\int_{\mathbb R^3}\partial_tu(x,t)\cdot u(x,t)dx\\
=2\int_{\mathbb R^3}(T-t)^{-\frac32}\{\partial_\tau
w(\frac{x}{(T-t)^{1/2}},\tau)+\frac{x}{2(T-t)^{1/2}}\cdot\nabla_yw(\frac{x}{(T-t)^{1/2}},\tau)\\
+\frac12w(\frac{x}{(T-t)^{1/2}},\tau)\}\cdot(T-t)^{-\frac12}w(\frac{x}{(T-t)^{1/2}},\tau)dx\\
=(T-t)^{-\frac12}\int_{\mathbb R^3}\{2\partial_\tau w(y,\tau)\cdot
w(y,\tau)+(y\cdot\nabla_y w(y,\tau))\cdot w(y,\tau)+|w(y,\tau)|^2dy.
\end{aligned}\end{equation} By using (2.13), from (2.3) we get (2.5) again.
\end{remark}

\vspace{0.5cm}

\section{$(L^\infty,L^2)$-decomposition of $w$}

In this section we shall prove that $w$ can be decomposed as the
sum of a $L^\infty(0,\infty; L^m(\mathbb R^3))$ ($m\in
[4,\infty]$) part and a $L^\infty(0,\infty; L^2(\mathbb R^3))\cap
L^2(0,\infty; H^1(\mathbb R^3))$ part.\\

Let $\varphi\in C_0^\infty(\mathbb R^3,[0,1])$ be a radial
symmetrical function satisfying
\begin{equation}\label{3.1}
\varphi(\xi)=1\quad \forall |\xi|\leq 1,\quad
\varphi(\xi)=0\quad\forall |\xi|\geq 2,\quad
\xi\cdot\nabla\varphi(\xi)\leq 0\quad \forall \xi.
\end{equation}
Like the Littlewood-Paley analysis, we define the operators
$$
\Delta_{-1}f=\mathcal{F}^{-1}[\varphi(\xi)\mathcal{F}[f](\xi)],\quad
\Delta_0f=\mathcal{F}^{-1}[(1-\varphi(\xi))\mathcal{F}[f](\xi)].
$$

Denote
\begin{equation}\label{3.2}
\begin{aligned}
\underline{w}(y,\tau)=\Delta_{-1}w(y,\tau)=\mathcal{F}^{-1}[\varphi]\ast w(y,\tau),\\
\overline{w}(y,\tau)=w(y,\tau)-\underline{w}(y,\tau)=\Delta_0w(y,\tau)=\mathcal{F}^{-1}[1-\varphi]\ast
w (y,\tau),\\
\tilde{\overline{w}}(y,\tau)=\mathcal{F}^{-1}[\sqrt{1-\varphi^2}]\ast
w(y,\tau).
\end{aligned}
\end{equation}
 Notice that
$$
\|w(\tau)\|^2_{L^2(\mathbb
R^3)}=\|\underline{w}(\tau)\|_{L^2(\mathbb
R^3)}^2+\|\tilde{\overline{w}}(\tau)\|_{L^2(\mathbb R^3)}^2.
$$
So (\ref{2.5}) can be written as
\begin{equation}\label{3.3}
\begin{aligned}
\frac12\frac{d}{d\tau}\int_{\mathbb
R^3}|\tilde{\overline{w}}(y,\tau)|^2dy \leq -\int_{\mathbb
R^3}|\nabla
\tilde{\overline{w}}(y,\tau)|^2-\frac14|\tilde{\overline{w}}(y,\tau)|^2dy\\
-\int_{\mathbb R^3}|\nabla
\underline{w}(y,\tau)|^2-\frac14|\underline{w}(y,\tau)|^2dy\\
-\frac12\frac{d}{d\tau}\int_{\mathbb
R^3}|\underline{w}(y,\tau)|^2dy.
\end{aligned}
\end{equation}

Applying the operator $\Delta_{-1}$ to the first equation of
(\ref{2.2}), we have
\begin{equation}\label{3.4}
\partial_\tau\Delta_{-1}w=\Delta
\Delta_{-1}w-\frac12\Delta_{-1}(y\cdot\nabla
w)-\frac12\Delta_{-1}w-\Delta_{-1}((w\cdot\nabla)w)-\nabla\Delta_{-1}q.
\end{equation}
Multiplying (\ref{3.4}) by $\Delta_{-1}w$ and integrating over
$\mathbb R^3$ we get
\begin{equation}\label{3.5}
\begin{aligned}
\frac12\frac{d}{d\tau}\int_{\mathbb R^3}|\Delta_{-1}w|^2dy
=-\int_{\mathbb R^3}|\nabla\Delta_{-1}w|^2dy-\frac12\int_{\mathbb
R^3}|\Delta_{-1}w|^2dy\\
-\frac12\int_{\mathbb R^3}\Delta_{-1}(y\cdot\nabla
w)\cdot\Delta_{-1}wdy-\int_{\mathbb
R^3}\Delta_{-1}((w\cdot\nabla)w)\cdot\Delta_{-1}wdy,
\end{aligned}
\end{equation}
 where div $w=0$ is used to cancel the term including $q$.

Because
\begin{eqnarray*}
&&\int_{\mathbb R^3}y\cdot\nabla|\Delta_{-1}w|^2dy=2\int_{\mathbb
R^3}y_j\Delta_{-1}w\cdot\partial_j\Delta_{-1}wdy\\
&&=2\int_{\mathbb
R^3}\xi_j\varphi\mathcal{F}[w]\cdot\partial_j\overline{(\varphi\mathcal{F}[w])}d\xi\\
&&=-3\int_{\mathbb R^3}\varphi^2|\mathcal{F}[w]|^2dy,
\end{eqnarray*}
we have
$$
\int_{\mathbb R^3}\partial_j\{y_j|\Delta_{-1}w|^2\}dy=0.
$$
So
\begin{eqnarray*}
&&\int_{\mathbb R^3}\partial_j(\Delta_{-1}(y_jw)\cdot\Delta_{-1}w)dy\\
&&=\int_{\mathbb
R^3}\partial_j\{\mathcal{F}^{-1}[\varphi]\ast(y_jw)\cdot\mathcal{F}^{-1}[\varphi]\ast
w\}dy\\
&&=(-1)\int_{\mathbb R^3}\partial_j\{\tilde{\varphi}_j\ast
w\cdot\mathcal{F}^{-1}[\varphi]\ast
w\}dy+\int_{\mathbb R^3}\partial_j\{y_j|\Delta_{-1}w|^2\}dy\\
&&=0
\end{eqnarray*}
where $\tilde{\varphi}_j(y)=y_j\mathcal{F}^{-1}[\varphi](y)$.\\

Noting that
\begin{eqnarray*}
&&\int\Delta_{-1}(y\cdot\nabla w)\cdot\Delta_{-1}wdy\\
&&=-\sum_{j=1}^3\int\Delta_{-1}(y_jw)\cdot\Delta_{-1}\partial_jw
dy-3\int|\Delta_{-1}w|^2dy\\
&&=-\sum_{j=1}^3\int\varphi(\xi)\mathcal{F}[y_jw]\cdot
\overline{\varphi(\xi)\mathcal{F}[\partial_jw]}d\xi-3\int|\Delta_{-1}w|^2dy
\end{eqnarray*}
and
$$
\mathcal{F}[y_jw]=i\frac{\partial}{\partial\xi_j}\mathcal{F}[w],\quad
\mathcal{F}[\partial_jw]=i\xi_j\mathcal{F}[w],
$$
we have
\begin{equation}\label{3.6}
\begin{aligned}
\int\Delta_{-1}(y\cdot\nabla w)\cdot\Delta_{-1}wdy
=-\sum_{j=1}^3\int\varphi^2(\xi)\xi_j\frac{\partial}{\partial\xi_j}\mathcal{F}[w]\cdot\overline{\mathcal{F}[w]}
d\xi-3\int|\Delta_{-1}w|^2dy\\
=-\sum_{j=1}^3\frac12\int\varphi^2(\xi)\xi_j\frac{\partial}{\partial\xi_j}|\mathcal{F}[w]|^2d\xi-3\int|\Delta_{-1}w|^2dy\\
=\sum_{j=1}^3\frac12\int
\xi_j\frac{\partial}{\partial\xi_j}\varphi^2(\xi)|\mathcal{F}[w]|^2d\xi+\frac32\int
\varphi^2(\xi)|\mathcal{F}[w]|^2d\xi-3\int|\Delta_{-1}w|^2dy\\
=\frac12\int
\xi\cdot\nabla\varphi^2(\xi)|\mathcal{F}[w]|^2d\xi-\frac32\int|\Delta_{-1}w|^2dy.
\end{aligned}
\end{equation}
 From (\ref{3.5})-(\ref{3.6}), we get
\begin{equation}\label{3.7}
\begin{aligned}
\frac12\frac{d}{d\tau}\int_{\mathbb R^3}|\Delta_{-1}w|^2dy
=-\int_{\mathbb R^3}|\nabla\Delta_{-1}w|^2dy+\frac14\int_{\mathbb
R^3}|\Delta_{-1}w|^2dy\\
-\frac14\int_{\mathbb
R^3}\xi\cdot\nabla\varphi^2(\xi)|\mathcal{F}[w](\xi,\tau)|^2d\xi-\int_{\mathbb
R^3}\Delta_{-1}((w\cdot\nabla)w)\cdot\Delta_{-1}wdy.
\end{aligned}
\end{equation}
 From (\ref{3.3}) and (\ref{3.7}), we have
\begin{equation}\label{3.8}
\begin{aligned}
\frac12\frac{d}{d\tau}\int_{\mathbb
R^3}|\tilde{\overline{w}}(y,\tau)|^2dy \leq -\int_{\mathbb
R^3}|\nabla
\tilde{\overline{w}}(y,\tau)|^2-\frac14|\tilde{\overline{w}}(y,\tau)|^2dy\\
+\frac14\int_{\mathbb
R^3}\xi\cdot\nabla\varphi^2(\xi)|\mathcal{F}[w](\xi,\tau)|^2d\xi+\int_{\mathbb
R^3}\Delta_{-1}((w\cdot\nabla)w)\cdot\Delta_{-1}wdy\\
\leq-\frac34\int_{\mathbb R^3}|\nabla
\tilde{\overline{w}}(y,\tau)|^2dy\\
+\frac14\int_{\mathbb
R^3}\xi\cdot\nabla\varphi^2(\xi)|\mathcal{F}[w](\xi,\tau)|^2d\xi+\int_{\mathbb
R^3}\Delta_{-1}((w\cdot\nabla)w)\cdot\Delta_{-1}wdy
\end{aligned}
\end{equation}
 where $|\xi||\mathcal{F}[\tilde{\overline{w}}]|^2\geq
|\mathcal{F}[\tilde{\overline{w}}]|^2$ is used in the last step.\\

Let $\alpha\in (0,\frac18)$ and define
$$
\chi(\xi)=\left\{\aligned & |\xi|^{\frac12+2\alpha}\varphi(\xi),\quad \forall |\xi|\leq \frac12+\alpha\\
                          & (\frac12+\alpha)^{\frac12+2\alpha}\varphi(\xi),\quad \forall
                          |\xi|\geq  \frac12+\alpha.
                          \endaligned\right.
$$
Instead of $\varphi$ by $\chi$, we define the operator
$$
\tilde{\Delta}_{-1}f=\mathcal{F}^{-1}[\chi(\xi)\mathcal{F}[f](\xi)].
$$

Applying $\tilde{\Delta}_{-1}$ to (\ref{2.2}), as (\ref{3.7}) we
have
\begin{eqnarray*}
&&\frac12\frac{d}{d\tau}\int_{\mathbb R^3}|\tilde{\Delta}_{-1}w|^2dy\\
&&=-\int_{\mathbb
R^3}|\nabla\tilde{\Delta}_{-1}w|^2dy+\frac14\int_{\mathbb
R^3}|\tilde{\Delta}_{-1}w|^2dy\\
&&-\frac14\int_{\mathbb
R^3}\xi\cdot\nabla\chi^2(\xi)|\mathcal{F}[w](\xi,\tau)|^2d\xi-\int_{\mathbb
R^3}\tilde{\Delta}_{-1}((w\cdot\nabla)w)\cdot\tilde{\Delta}_{-1}wdy.
\end{eqnarray*}
Combining it with (\ref{3.8}), we have
\begin{eqnarray*}
&&\frac12\frac{d}{d\tau}\int_{\mathbb
R^3}|\tilde{\Delta}_{-1}w|^2+|\tilde{\overline{w}}|^2dy\\
&&\leq-\frac34\int_{\mathbb
R^3}|\nabla\tilde{\overline{w}}|^2dy-\alpha\int_{\mathbb
R^3}|\tilde{\Delta}_{-1}w|^2dy\\
&&+\int_{\mathbb
R^3}\Delta_{-1}((w\cdot\nabla)w)\cdot\Delta_{-1}w-\tilde{\Delta}_{-1}((w\cdot\nabla)w)\cdot\tilde{\Delta}_{-1}wdy\\
&&+\frac14\int_{\mathbb
R^3}\xi\cdot\nabla(\varphi^2-\chi^2)|\mathcal{F}[w]|^2d\xi-\int_{\mathbb
R^3}|\nabla\tilde{\Delta}_{-1}w|^2dy
+(\frac14+\alpha)\int_{\mathbb R^3}|\tilde{\Delta}_{-1}w|^2dy.
\end{eqnarray*}
For $|\xi|\leq 1$, the last term is written as
$$
A=\int\{\frac14\xi\cdot\nabla(-\chi^2)-|\xi|^2\chi^2+(\frac14+\alpha)\chi^2\}|\mathcal{F}[w]|^2d\xi
$$
and noting that $\varphi(\xi)=1$ for $|\xi|\leq1$ as well as the
definition of $\chi$,  $A\leq 0$. For $|\xi|\in[1,2]$, the last
term is written as
$$
B=\int\{\frac14(1-(\frac12+\alpha)^{1+4\alpha})\xi\cdot\nabla\varphi^2
-(|\xi|^2-(\frac14+\alpha))(\frac12+\alpha)^{1+4\alpha}\varphi^2\}|\mathcal{F}[w]|^2d\xi,
$$
and $B\leq 0$. So we get
\begin{equation}\label{3.9}
\begin{aligned}
\frac12\frac{d}{d\tau}\int_{\mathbb
R^3}|\tilde{\Delta}_{-1}w|^2+|\tilde{\overline{w}}|^2dy
\leq-\frac34\int_{\mathbb
R^3}|\nabla\tilde{\overline{w}}|^2dy-\alpha\int_{\mathbb
R^3}|\tilde{\Delta}_{-1}w|^2dy\\
+\int_{\mathbb
R^3}\Delta_{-1}((w\cdot\nabla)w)\cdot\Delta_{-1}w-\tilde{\Delta}_{-1}((w\cdot\nabla)w)\cdot\tilde{\Delta}_{-1}wdy.
\end{aligned}
\end{equation}

\begin{lemma}\label{lem3.1} (1) For any $m\in[4,\infty]$,
$$
\|\Delta_{-1}f\|_{L^m(\mathbb R^3)}\leq C(\alpha)
\|\tilde{\Delta}_{-1}f\|_{L^2(\mathbb R^3)},\quad\forall f\in
L^2(\mathbb R^3)
$$
where the constant $C(\alpha)<\infty$ depends only on $\alpha$.

(2) For all $\beta=(\beta_1,\beta_2,\beta_3)$ ($\beta_j\in\mathbb
N$, $j=1,2,3$)
$$
\|D^\beta\Delta_0w(\cdot,\tau)\|_{L^2(\mathbb R^3)}\leq
\|D^\beta\tilde{\overline{w}}(\cdot,\tau)\|_{L^2(\mathbb R^3)}
$$
where
$D^\beta=\partial_1^{\beta_1}\partial_2^{\beta_2}\partial_3^{\beta_3}$.

\end{lemma}

\vspace{0.5cm}

$Proof.$ From  Hausdorff-Young inequality
\begin{eqnarray*}
&&\|\Delta_{-1}f\|_{L^m(\mathbb R^3)}\\
&&\leq
(2\pi)^{3/m'}(\int_{\mathbb R^3}|\varphi(\xi)\mathcal{F}[f](\xi)|^{m'}d\xi)^{1/m'}\\
&&\leq (2\pi)^{3/m'}(\int_{\mathbb
R^3}||\xi|^{\frac12+2\alpha}\varphi(\xi)\mathcal{F}[f](\xi)|^2)^{1/2}(\int_{|\xi|\leq
2}|\xi|^{-(\frac12+2\alpha)\frac{2m'}{2-m'}}d\xi)^{\frac{2-m'}{2m'}}\\
&&\leq C(\alpha)(\int_{\mathbb
R^3}|\chi(\xi)\mathcal{F}[f](\xi)|^2)^{1/2}
\end{eqnarray*}
for $\alpha\in(0,\frac18)$. So we have (1).

To prove (2), we only need to consider the case $|\beta|=\sum_{1\leq
j\leq 3}\beta_j=0$. Since $0\leq\varphi\leq 1$ and
$1-\varphi^2=(1-\varphi)(1+\varphi)\geq (1-\varphi)^2$, in this case
we have
\begin{eqnarray*}
&&\int_{\mathbb R^3}|\Delta_0w(y,\tau)|^2dy=\int_{\mathbb
R^3}(1-\varphi(\xi))^2|\mathcal{F}[w](\xi,\tau)|^2d\xi\\
&&\leq\int_{\mathbb
R^3}(1-\varphi^2(\xi))|\mathcal{F}[w](\xi,\tau)|^2d\xi=\int_{\mathbb
R^3}|\tilde{\overline{w}}(y,\tau)|^2dy.\qed
\end{eqnarray*}

\vspace{0.5cm}

Now we estimate the last term in the right of (\ref{3.9}). We only
need to consider the integration for the first function in the
last term, because for another function the proof is same.  Notice
that
\begin{eqnarray*}
&&\int \Delta_{-1}(w_j\partial_jw)\cdot\Delta_{-1}w
dy=-\int\Delta_{-1}(w_jw)\cdot\Delta_{-1}(\partial_jw)dy\\
&&=-\int\Delta_{-1}(\Delta_{-1}w_j
\Delta_{-1}w)\cdot\Delta_{-1}(\partial_jw)dy-\int\Delta_{-1}(\Delta_{-1}w_j\Delta_0w)\cdot\Delta_{-1}\partial_jwdy\\
&&-\int\Delta_{-1}(\Delta_0w_j\Delta_{-1}w)\cdot\Delta_{-1}(\partial_jw)dy-\int\Delta_{-1}(\Delta_0w_j\Delta_0w)\cdot\Delta_{-1}(\partial_jw)dy.
\end{eqnarray*}
Because
\begin{eqnarray*}
&&|\int\Delta_{-1}(\Delta_{-1}w_j
\Delta_{-1}w)\cdot\Delta_{-1}(\partial_jw)dy|\\
&&\leq(\int|\Delta_{-1}w|^4dy)^{1/2}(\int\varphi^2(\xi)|\xi_j\mathcal{F}[w](\xi,\tau)|^2dx)^{1/2}\\
&&\leq C\|\tilde{\Delta}_{-1}w\|_{L^2(\mathbb R^3)}^3,
\quad\text{(by Lemma \ref{lem3.1} (1) and the definition of
$\tilde{\Delta}_{-1}$)}
\end{eqnarray*}
and
\begin{eqnarray*}
&&|\int\Delta_{-1}(\Delta_{-1}w_j\Delta_0w)\cdot\Delta_{-1}\partial_jwdy|\\
&&\leq (\int|\Delta_{-1}w|^4dy)^{1/2}(\int|\Delta_0w|^2dy)^{1/2}\\
&&\leq C\|\tilde{\Delta}_{-1}w\|_{L^2(\mathbb
R^3)}^2\|\tilde{\overline{w}}\|_{L^2(\mathbb R^3)}\quad\text{(by
Lemma \ref{lem3.1} (1)-(2))}
\end{eqnarray*}
as well as
\begin{eqnarray*}
&&|\int\Delta_{-1}(\Delta_0w_j\Delta_0w)\cdot\Delta_{-1}(\partial_jw)dy|\\
&&\leq \|\Delta_{-1}(\partial_jw)\|_{L^\infty(\mathbb
R^3)}\|\Delta_0w\|_{L^2(\mathbb R^3)}^2\\
&&\leq C\|\tilde{\Delta}_{-1}w\|_{L^2(\mathbb
R^3)}\|\tilde{\overline{w}}\|_{L^2(\mathbb R^3)}^2\quad\text{(by
Lemma \ref{lem3.1} (1)-(2))}
\end{eqnarray*}
the last term in the right of (\ref{3.9}) can be estimated by
$$
C\{\|\tilde{\Delta}_{-1}w\|_{L^2(\mathbb R^3)}^3+
\|\tilde{\Delta}_{-1}w\|_{L^2(\mathbb
R^3)}^2\|\tilde{\overline{w}}\|_{L^2(\mathbb R^3)}+
\|\tilde{\Delta}_{-1}w\|_{L^2(\mathbb
R^3)}\|\tilde{\overline{w}}\|_{L^2(\mathbb R^3)}^2 \}.
$$
So we get
\begin{equation}\label{3.10}
\begin{aligned}
\frac12\frac{d}{d\tau}\int_{\mathbb
R^3}|\tilde{\Delta}_{-1}w(y,\tau)|^2+|\tilde{\overline{w}}(y,\tau)|^2dy+(\frac34-\alpha)\int_{\mathbb R^3}|\nabla\tilde{\overline{w}}(y,\tau)|^2dy\\
\leq-\alpha\int_{\mathbb
R^3}|\tilde{\Delta}_{-1}w(y,\tau)|^2+|\tilde{\overline{w}}(y,\tau)|^2dy\\
+C(\int_{\mathbb
R^3}|\tilde{\Delta}_{-1}w(y,\tau)|^2+|\tilde{\overline{w}}(y,\tau)|^2dy)^{3/2}
\end{aligned}
\end{equation}

\begin{proposition}\label{pro3.2}There is $\delta>0$ such that if
\begin{equation}\label{3.11}
\int_{\mathbb
R^3}|\tilde{\Delta}_{-1}w(y,0)|^2+|\tilde{\overline{w}}(y,0)|^2dy\leq\delta
\end{equation}
 then for all $\tau>0$
\begin{equation}\label{3.12}
\frac{d}{d\tau}\int_{\mathbb
R^3}|\tilde{\Delta}_{-1}w(y,\tau)|^2+|\tilde{\overline{w}}(y,\tau)|^2dy\leq-\alpha\int_{\mathbb
R^3}|\tilde{\Delta}_{-1}w(y,\tau)|^2+|\tilde{\overline{w}}(y,\tau)|^2dy.
\end{equation}
 Moreover
$w(y,\tau)=\underline{w}(y,\tau)+\overline{w}(y,\tau)$, and for
all $m\in[4,\infty]$,
\begin{equation}\label{3.13}
\begin{aligned}
\|D^\beta\underline{w}(\tau)\|_{L^m(\mathbb R^3)}\leq
C(\beta)\delta,\quad\forall\tau>0,\quad \forall \beta,\\
\lim_{\tau\to\infty}\|\underline{w}(\tau)\|_{L^m(\mathbb R^3)}=0,
\end{aligned}
\end{equation}
\begin{equation}\label{3.14}
\begin{aligned}
\sup_{\tau\geq0}\int_{\mathbb
R^3}|\overline{w}(y,\tau)|^2dy+\int_0^\infty d\tau\int_{\mathbb
R^3}|\nabla\overline{w}(y,\tau)|^2dy\leq C\delta,\\
\lim_{\tau\to\infty}\int_{\mathbb
R^3}|\overline{w}(y,\tau)|^2dy=0.
\end{aligned}
\end{equation}

\end{proposition}

\vspace{0.5cm}

For example, we may take $\delta\leq (\frac{\alpha}{2C})^2$.
Proposition \ref{pro3.2} follows from (\ref{3.10}) and Lemma
\ref{lem3.1}. Note that
\begin{equation}\label{3.15}
\begin{aligned}
\int_{\mathbb
R^3}|\tilde{\Delta}_{-1}w(y,0)|^2+|\tilde{\overline{w}}(y,0)|^2dy
=\int_{\mathbb
R^3}(\chi^2(\xi)+1-\varphi^2(\xi))|\mathcal{F}[w](\xi,0)|^2d\xi\\
\leq \int_{\mathbb R^3}|\mathcal{F}[w](\xi,0)|^2d\xi=\int_{\mathbb
R^3}|w(y,0)|^2dy=\int_{\mathbb R^3}|u_0(x)|^2dx.
\end{aligned}
\end{equation}
 So we have

\begin{corollary}\label{cor3.3} There is $\delta>0$ such that if
$\|u_0\|_{L^2(\mathbb R^3)}\leq\delta^{1/2}$, then we have
(\ref{3.12})-(\ref{3.14}).
\end{corollary}

\vspace{0.5cm}

\begin{remark}\label{rmk3.4} Suppose $\psi$ is a function satisfying
\begin{equation}\label{3.16}
\psi\in C(\mathbb R^3,[0,1]),\quad \xi\cdot\nabla_\xi\psi(\xi)\in
L^\infty(\mathbb R^3).
\end{equation}
Since $\psi(\xi)\mathcal{F}[w](\xi,\tau)\in L^2(\mathbb R^3)$, we
have $\mathcal{F}^{-1}[\psi]\ast
w=\mathcal{F}^{-1}[\psi\mathcal{F}[w]]\in L^2(\mathbb R^3)$ and
\begin{equation}\label{3.17}
\begin{aligned}
\int_{\mathbb R^3}|\mathcal{F}^{-1}[\psi]\ast w(y,\tau)|^2dy(T-t)^{\frac32}\\
=\int_{\mathbb R^3}|\mathcal{F}^{-1}[\psi]\ast
w(\frac{\mu}{(T-t)^{1/2}},\tau)|^2d\mu\\
=\int_{\mathbb R^3}|\int_{\mathbb
R^3}\mathcal{F}^{-1}[\psi](\frac{\mu}{(T-t)^{1/2}}-z)w(z,\tau)dz|^2d\mu\\
=\int_{\mathbb R^3}|\int_{\mathbb
R^3}\mathcal{F}[\psi](\frac{\mu}{(T-t)^{1/2}}-z)(T-t)^{1/2}u((T-t)^{1/2}z,t)dz|^2d\mu\\
=\int_{\mathbb R^3}|\int_{\mathbb
R^3}\mathcal{F}^{-1}[\psi](\frac{\mu-x}{(T-t)^{1/2}})u(x,t)dx|^2d\mu
(T-t)^{-2}
\end{aligned}\end{equation}
where $\tau=(-)\ln(T-t)$. Note that
\begin{equation}\label{3.18}
\begin{aligned}
\partial_t\{(T-t)^{3/2}\psi((T-t)^{1/2}\xi)\mathcal{F}[u](\xi,t)\}\\
=\partial_t\mathcal{F}[\int_{\mathbb
R^3}\mathcal{F}^{-1}[\psi](\frac{\mu-x}{(T-t)^{1/2}})u(x,t)dx]\\
=\mathcal{F}[\int_{\mathbb
R^3}\mathcal{F}^{-1}[\psi](\frac{\mu-x}{(T-t)^{1/2}})\partial_tu(x,t)+\{\frac{\mu-x}{2(T-t)^{3/2}}\cdot\mathcal{F}^{-1}[\psi]'(\frac{\mu-x}{(T-t)^{1/2}})\}u(x,t)dx]\\
=\mathcal{F}[\int_{\mathbb
R^3}\mathcal{F}^{-1}[\psi](\frac{\mu-x}{(T-t)^{1/2}})\{\partial_\tau
w(\frac{x}{(T-t)^{1/2}},\tau)
+\frac{x}{2(T-t)^{1/2}}\cdot\nabla_yw(\frac{x}{(T-t)^{1/2}},\tau)\\
+\frac12 w(\frac{x}{(T-t)^{1/2}},\tau) \}(T-t)^{-\frac32}\\
+\{\frac{\mu-x}{2(T-t)^{3/2}}\cdot\mathcal{F}^{-1}[\psi]'(\frac{\mu-x}{(T-t)^{1/2}})\}w(\frac{x}{(T-t)^{1/2}},\tau)(T-t)^{-\frac12}dx]
\end{aligned}\end{equation}
and
\begin{equation}\label{3.19}
\begin{aligned}
\frac{d}{dt}\int_{\mathbb R^3}|\int_{\mathbb
R^3}\mathcal{F}^{-1}[\psi](\frac{\mu-x}{(T-t)^{1/2}})u(x,t)dx|^2d\mu\\
=\frac{d}{dt}\int_{\mathbb
R^3}|(T-t)^{3/2}\psi((T-t)^{1/2}\xi)\mathcal{F}[u](\xi,t)|^2d\xi\\
=\frac{d}{dt}\int_{\mathbb
R^3}\int_0^t\partial_h|(T-h)^{3/2}\psi((T-h)^{1/2}\xi)\mathcal{F}[u](\xi,h)|^2dhd\xi\\
=\frac{d}{dt}\int_0^t\int_{\mathbb
R^3}\partial_h|(T-h)^{3/2}\psi((T-h)^{1/2}\xi)\mathcal{F}[u](\xi,h)|^2d\xi dh\\
=\int_{\mathbb
R^3}\partial_t|(T-t)^{3/2}\psi((T-t)^{1/2}\xi)\mathcal{F}[u](\xi,t)|^2d\xi\\
=2\int_{\mathbb
R^3}\{(-)\frac32(T-t)^{1/2}\psi((T-t)^{1/2}\xi)\mathcal{F}[u](\xi,t)-(T-t)\frac{\xi}2\cdot\psi'((T-t)^{1/2}\xi)\mathcal{F}[u](\xi,t)\\
+(T-t)^{3/2}\psi((T-t)^{1/2}\xi)\partial_t\mathcal{F}[u](\xi,t)\}\cdot
(T-t)^{3/2}\psi((T-t)^{1/2}\xi)\overline{\mathcal{F}[u]}(\xi,t)d\xi
\end{aligned}\end{equation}
where noting that $\mathcal{F}[u](\xi,t),\,\,
\partial_t\mathcal{F}[u](\xi,t)\in L^2(\mathbb R^3)$ for $t<T$ and
$\psi$ satisfies (3.16), we have
$$
\int_{\mathbb
R^3}|\partial_t|(T-t)^{3/2}\psi((T-t)^{1/2}\xi)\mathcal{F}[u](\xi,t)|^2|d\xi<\infty
$$
and Fubini theorem can be used.\\

From (\ref{3.17})-(\ref{3.19}), we get
\begin{equation}\label{3.20}
\begin{aligned}
(T-t)^{5/2}\frac{d}{d\tau}\int_{\mathbb
R^3}|\mathcal{F}^{-1}[\psi]\ast
w(y,\tau)-\frac{7}{2}(T-t)^{5/2}\int_{\mathbb
R^3}|\mathcal{F}^{-1}[\psi]\ast w(y,\tau)|^2dy\\
=\frac{d}{dt}\{(T-t)^{7/2}\int_{\mathbb
R^3}|\mathcal{F}^{-1}[\psi]\ast w(y,\tau)|^2dy\}\\
=\frac{d}{dt}\int_{\mathbb R^3}|\int_{\mathbb
R^3}\mathcal{F}^{-1}[\psi](\frac{\mu-x}{(T-t)^{1/2}})u(x,t)dx|^2d\mu\\
=2(T-t)^{5/2}\int_{\mathbb
R^3}\{(-)\frac32\psi(\xi)\mathcal{F}[w](\xi,\tau)-\frac{\xi}2\cdot\psi'(\xi)\mathcal{F}[w](\xi,\tau)\}\cdot\psi(\xi)\overline{\mathcal{F}[w]}(\xi,\tau)d\xi  \\
+2\int_{\mathbb R^3}\{\int_{\mathbb
R^3}\mathcal{F}^{-1}[\psi](\frac{\mu-x}{(T-t)^{1/2}})\partial_tu(x,t)dx\}\cdot\{\int_{\mathbb
R^3}\mathcal{F}^{-1}[\psi](\frac{\mu-x}{(T-t)^{1/2}})u(x,t)dx\}d\mu\\
=(-3)(T-t)^{5/2}\int_{\mathbb R^3}|\mathcal{F}^{-1}[\psi]\ast w(y,\tau)|^2dy\\
-\frac{(T-t)^{5/2}}2\int_{\mathbb
R^3}(\xi\cdot\nabla_\xi\psi^2(\xi))|\mathcal{F}[w](\xi,\tau)|^2d\xi\\
+2(T-t)^{5/2}\int_{\mathbb R^3}\int_{\mathbb
R^3}\mathcal{F}^{-1}[\psi](y-z)\{\partial_\tau w(z,\tau)
+\frac{z}{2}\cdot\nabla_zw(z,\tau) +\frac12 w(z,\tau) \}dz\\
\cdot\{\mathcal{F}^{-1}[\psi]\ast w(y,\tau)\}dy
\end{aligned}\end{equation}
So we have
\begin{equation}\label{3.21}
\begin{aligned}
2\int_{\mathbb R^3}\mathcal{F}^{-1}[\psi]\ast\{\partial_\tau
w+\frac{y}2\cdot\nabla_yw\}(y,\tau)\cdot\mathcal{F}^{-1}[\psi]\ast
w(y,\tau)dy\\
=\frac{d}{d\tau}\int_{\mathbb R^3}|\mathcal{F}^{-1}[\psi]\ast
w(y,\tau)|^2dy-\frac32\int_{\mathbb R^3}|\mathcal{F}^{-1}[\psi]\ast
w(y,\tau)|^2dy\\
+\frac12\int_{\mathbb
R^3}(\xi\cdot\nabla_\xi\psi^2(\xi))|\mathcal{F}[w](\xi,\tau)|^2d\xi.
\end{aligned}\end{equation}

Note that $\varphi$ and $\chi$ satisfy (3.16), and we can use (3.21)
to obtain (3.7) for $\varphi$ and $\chi$ again. Furthermore, notice
that $1-\varphi$ satisfies (3.16) and
$\|\partial_t\nabla_xu(t)\|_{L^2(\mathbb R^3)}$ is bounded for
$t<T$, we can prove  the same equation as (3.21) for $(1-\varphi)$
and $\nabla_yw$ instead of $\psi$ and $w$, which can be used to
obtain (4.4) of section 4 from (4.1) too.

\end{remark}

\vspace{0.5cm}

\section {$L^\infty$-estimate of $\overline{w}$}

Applying the operator $\Delta_0$ (see (\ref{3.2})) to (\ref{2.6}),
and integrating over $\mathbb R^3$ we have
\begin{equation}\label{4.1}
\begin{aligned}
\frac12\frac{d}{d\tau}\int_{\mathbb R^3}|\Delta_0\nabla
w|^2dy=-\int_{\mathbb R^3}|\nabla^2\Delta_0w|^2dy-\int_{\mathbb
R^3}|\Delta_0\nabla w|^2dy\\
-\sum_{j=1}^3\frac12\int_{\mathbb
R^3}\Delta_0(y\cdot\nabla\partial_jw)\cdot\Delta_0\partial_jw dy\\
-\sum_{j=1}^3\int_{\mathbb
R^3}\Delta_0((\partial_jw\cdot\nabla)w)\cdot\Delta_0\partial_jw+\Delta_0((w\cdot\nabla)\partial_jw)\cdot\Delta_0\partial_jw
dy.
\end{aligned}
\end{equation}
 Since the support set of $1-\varphi$ is not compact, we can
not do the same thing as in (\ref{3.6}) for the 3rd term in the
right side of (\ref{4.1}). But with more patient, by using
$\Delta_0f=f-\Delta_{-1}f$, we have
\begin{equation}\label{4.2}
\begin{aligned}
\int\Delta_0((y\cdot\nabla)\partial_jw)\cdot\Delta_0\partial_jwdy\\
=\int((y\cdot\nabla)\partial_jw)\cdot\partial_jw
dy-\int((y\cdot\nabla)\partial_jw)\cdot\Delta_{-1}\partial_jwdy\\
-\int\Delta_{-1}((y\cdot\nabla)\partial_jw)\cdot\partial_jwdy+\int\Delta_{-1}((y\cdot\nabla)\partial_jw)\cdot\Delta_{-1}(\partial_jw)dy.
\end{aligned}
\end{equation}
 As in (\ref{2.4}), we have
$$
\int((y\cdot\nabla)\partial_jw)\cdot\partial_jw dy\geq-\frac32\int
|\partial_jw|^2dy.
$$
On the other hand, as in (\ref{3.6}) we have
$$
\int\Delta_{-1}((y\cdot\nabla)\partial_jw)\cdot\Delta_{-1}(\partial_jw)dy
=\frac12\int\xi\cdot\nabla\varphi^2|\mathcal{F}[\partial_jw]|^2d\xi-\frac32\int\varphi^2|\mathcal{F}[\partial_jw]|^2d\xi.
$$
The remainder in the right of (\ref{4.2}) is
\begin{eqnarray*}
&&2\int\varphi\partial_k(\xi_k\mathcal{F}[\partial_jw])\cdot\overline{\mathcal{F}[\partial_jw]}d\xi\\
&&=-\int2(\xi\cdot\nabla\varphi)|\mathcal{F}[\partial_jw]|^2+\varphi\xi\cdot\nabla|\mathcal{F}[\partial_jw]|^2d\xi\\
&&=-\int(\xi\cdot\nabla\varphi)|\mathcal{F}[\partial_jw]|^2d\xi+3\int\varphi|\mathcal{F}[\partial_jw]|^2d\xi.
\end{eqnarray*}
Then the right of (\ref{4.2}) is larger than
\begin{equation}\label{4.3}
\begin{aligned}
-\frac32\int
|\partial_jw|^2dy-\frac32\int\varphi^2|\mathcal{F}[\partial_jw]|^2d\xi
+\frac12\int\xi\cdot\nabla\varphi^2|\mathcal{F}[\partial_jw]|^2d\xi\\
-\int(\xi\cdot\nabla\varphi)|\mathcal{F}[\partial_jw]|^2d\xi+3\int\varphi|\mathcal{F}[\partial_jw]|^2d\xi\\
=\frac12\int\xi\cdot\nabla(1-\varphi(\xi))^2|\mathcal{F}[\partial_jw]|^2d\xi
-\frac32\int|\Delta_0\partial_jw|^2dy.
\end{aligned}
\end{equation}
Since from (\ref{3.1})
$$
\xi\cdot\nabla(1-\varphi(\xi))^2=|\xi|\frac{d}{d|\xi|}(1-\varphi(\xi))^2\geq
0 $$ Instead of the 3rd term in the right side of (\ref{4.1}) by
(\ref{4.2})-(\ref{4.3}), we get
\begin{equation}\label{4.4}
\begin{aligned}
\frac12\frac{d}{d\tau}\int_{\mathbb R^3}|\Delta_0\nabla
w|^2dy\leq-\int_{\mathbb
R^3}|\nabla^2\Delta_0w|^2dy-\frac14\int_{\mathbb
R^3}|\Delta_0\nabla w|^2dy\\
-\sum_{j=1}^3\int_{\mathbb
R^3}\Delta_0((\partial_jw\cdot\nabla)w)\cdot\Delta_0\partial_jw+\Delta_0((w\cdot\nabla)\partial_jw)\cdot\Delta_0\partial_jw
dy.
\end{aligned}
\end{equation}

Decompose the last integration of the right side of (\ref{4.4}) by
$w=\underline{w}+\overline{w}$ and note that
\begin{eqnarray*}
|\int((\overline{w}\cdot\nabla)\partial_j\overline{w})\cdot\partial_j\overline{w}dy|
&\leq&\|\nabla^2\overline{w}\|_{L^2(\mathbb
R^3)}(\int|\overline{w}|^2|\nabla\overline{w}|^2dy)^{1/2}\\
&\leq& C\|\nabla^2\overline{w}\|_{L^2(\mathbb
R^3)}^{3/2}\|\nabla\overline{w}\|_{L^2(\mathbb R^3)}^{3/2},
\end{eqnarray*}
$$
|\int((\underline{w}\cdot\nabla)\partial_j\overline{w})\cdot\partial_j\overline{w}dy|
\leq\|\underline{w}\|_{L^\infty(\mathbb
R^3)}\|\nabla^2\overline{w}\|_{L^2(\mathbb
R^3)}\|\nabla\overline{w}\|_{L^2(\mathbb R^3)},
$$
$$
|\int((\overline{w}\cdot\nabla)\partial_j\underline{w})\cdot\partial_j\overline{w}dy|\leq
C\|\underline{w}\|_{L^\infty(\mathbb
R^3)}\|\overline{w}\|_{L^2(\mathbb
R^3)}\|\nabla\overline{w}\|_{L^2(\mathbb R^3)},
$$
and
$$
|\int((\underline{w}\cdot\nabla)\partial_j\underline{w})\cdot\partial_j\overline{w}dy|
\leq
C(\int|\underline{w}|^4dy)^{1/2}\|\nabla\overline{w}\|_{L^2(\mathbb
R^3)}
$$
as well as the same estimates for another one. Then by Proposition
\ref{pro3.2}, we have
\begin{equation}\label{4.5}
\begin{aligned}
\frac12\frac{d}{d\tau}\int_{\mathbb R^3}|\nabla\overline{
w}(y,\tau)|^2dy\leq-\frac12\int_{\mathbb
R^3}|\nabla^2\overline{w}(y,\tau)|^2dy-\frac18\int_{\mathbb
R^3}|\nabla\overline{ w}(y,\tau)|^2dy\\
+C\|\nabla\overline{w}(\tau)\|_{L^2(\mathbb
R^3)}\{C\delta-\|\nabla\overline{w}(\tau)\|_{L^2(\mathbb
R^3)}+C\|\nabla\overline{w}(\tau)\|_{L^2(\mathbb R^3)}^5\}
\end{aligned}
\end{equation}
 Note that (see Remark \ref{rmk4.4}) there is $\delta_1>0$ such that if for some
$\tau_0\geq0$
\begin{equation}\label{4.6}
\|\nabla\overline{w}(\tau_0)\|_{L^2(\mathbb R^3)}\leq\delta_1
\end{equation}
then
$$
\|\nabla\overline{w}(\tau)\|_{L^2(\mathbb
R^3)}\leq\delta_1,\quad\forall \tau\geq\tau_0.
$$
From (\ref{3.14}), (\ref{4.6}) can be satisfied provided that
(\ref{3.11}) is satisfied. So we have

\begin{lemma}\label{lem4.1} Suppose (\ref{3.11}) is satisfied. Then there is
$\delta_1>0$ ($\delta_1\downarrow0$ as $\delta\downarrow0$) and
$\tau_0>0$ such that
$$
\|\nabla\overline{w}(\tau)\|_{L^2(\mathbb
R^3)}\leq\delta_1,\quad\forall \tau\geq\tau_0.
$$
\end{lemma}

\vspace{0.5cm}

Estimate the last term in the right side of (\ref{2.7}) by using
$w=\underline{w}+\overline{w}$, and note that
\begin{eqnarray*}
|\int\partial_j\underline{w}_k\partial_j\underline{w}_l\partial_l\underline{w}_kdy|
&\leq& \|\nabla\underline{w}\|_{L^\infty(\mathbb R^3)}\int|\nabla\underline{w}|^2dy\\
&\leq& C\delta\int|\nabla\underline{w}|^2dy,
\end{eqnarray*}
\begin{eqnarray*}
|\int\partial_j\underline{w}_k\partial_j\underline{w}_l\partial_l\overline{w}_kdy|
&\leq&(\int|\nabla\underline{w}|^4dy)^{1/2}(\int|\nabla\overline{w}|^2dy)^{1/2}\\
&\leq& C\delta(\int|\nabla\overline{w}|^2dy)^{1/2},
\end{eqnarray*}
and
\begin{eqnarray*}
|\int\partial_j\underline{w}_k\partial_j\overline{w}_l\partial_l\overline{w}_kdy|
&\leq&\|\nabla\underline{w}\|_{L^\infty(\mathbb
R^3)}\int|\nabla\overline{w}|^2dy\\
&\leq& C\delta\int|\nabla\overline{w}|^2dy,
\end{eqnarray*}
as well as
\begin{eqnarray*}
|\int\partial_j\overline{w}_k\partial_j\overline{w}_l\partial_l\overline{w}_kdy|
&\leq& C\|\nabla\overline{w}\|_{L^2(\mathbb
R^3)}^{3/2}\|\nabla^2\overline{w}\|_{L^2(\mathbb R^3)}^{3/2}\\
&\leq& \|\nabla^2\overline{w}\|_{L^2(\mathbb R^3)}^2+C\delta_1.
\end{eqnarray*}
So we have
\begin{equation}\label{4.7}
\frac{d}{d\tau}\int_{\mathbb R^3}|\nabla w(y,\tau)|^2dy \leq
-\int_{\mathbb R^3}|\nabla^2w(y,\tau)|^2dy-\frac12\int_{\mathbb
R^3}|\nabla w(y,\tau)|^2dy+C\delta_1.
\end{equation}

\begin{lemma}\label{lem4.2}  Suppose (\ref{3.11}) is satisfied. Then there is
$\delta_1>0$ ($\delta_1\downarrow0$ as $\delta\downarrow0$) and
$\tau_0>0$ such that for all $\tau\geq\tau_0$,
\begin{equation}\label{4.8}
\int_{\mathbb R^3}|\nabla w(y,\tau)|^2dy \leq
e^{-\frac12(\tau-\tau_0)}\int_{\mathbb R^3}|\nabla
w(y,\tau_0)|^2dy+2C\delta_1(1-e^{-\frac12(\tau-\tau_0)}).
\end{equation}

\end{lemma}
\vspace{0.5cm}

Considering  (\ref{2.8}), and noting that div $\Delta w=0$ implies
$$
\int(\Delta w)\cdot ((w\cdot\nabla)\Delta w)dy=0,
$$
the last term in the right side of (\ref{2.8}) can be written as the
sum of the following terms
$$
\int |\nabla^2w|^2|\nabla w|dy.
$$
Since it can be estimated by
\begin{eqnarray*}
&&(\int|\nabla w|^2dy)^{1/2}(\int|\nabla^2w|^4dy)^{1/2}\\
&&\leq C(\int|\nabla w|^2dy)^{1/2}(\int|\Delta
w|^2dy)^{1/4}(\int|\nabla\Delta w|^2dy)^{3/4}
\end{eqnarray*}
by (\ref{2.8}) and Lemma \ref{lem4.2} we have

\begin{lemma}\label{lem4.3} Suppose (\ref{3.11}) is satisfied. Then there is
$\delta_1>0$ ($\delta_1\downarrow0$ as $\delta\downarrow0$) and
$\tau_0>0$ such that for all $\tau\geq\tau_0$,
$$
\int_{\mathbb R^3}|\Delta w(y,\tau)|^2dy\leq
e^{-(\frac32-C\delta_1)(\tau-\tau_0)}\int_{\mathbb R^3}|\Delta
w(y,\tau_0)|^2dy
$$

\end{lemma}

\vspace{0.5cm}

From Lemma \ref{lem4.1}-\ref{lem4.3} and Corollary \ref{cor3.3},
we proved the Theorem \ref{thm1.1}.

\begin{remark}\label{rmk4.4} Suppose a nonnegative continuous function
$h(\tau)$ satisfies
$$
\frac{d}{d\tau}h(\tau)\leq F(h(\tau)):=
C\delta-Bh(\tau)+h^5(\tau),\quad\forall\tau>0,
$$
where $C$, $B$ and $\delta$ are positive constants. If $\delta$ is
small enough so that
$$
h_-:=\frac12(B-\sqrt{B^2-4C\delta})\in (0,1),
$$
and if $h(0)< h_-$, then for all $\tau>0$, $F(h(\tau))\leq
C\delta-Bh(\tau)+h^2(\tau)$ and
$$
h(\tau)\in [0,h_-].
$$
\end{remark}

\end{document}